\documentclass[12pt,reqno]{amsart}
\usepackage{amsmath,amsfonts,amssymb,amsthm,amsxtra,amscd}
\begin{document}

\title[Stably extendable vector bundles]{Infinitely stably extendable 
       vector bundles on projective spaces}

\author[I. Coand\u{a}]{Iustin Coand\u{a}}
\address{Institute of Mathematics of the Romanian Academy, 
         P. O. Box 1-764, RO-014700, Bucharest, Romania}
\email{Iustin.Coanda@imar.ro}

\subjclass[2000]{Primary: 14F05; Secondary: 14J60, 14B10}

\keywords{Vector bundle, projective space, monad, extendable vector bundle}

\thanks{Partially supported by CNCSIS grant ID-PCE no.51/28.09.2007}

\begin{abstract}
According to Horrocks (1966), a vector bundle $E$ on the projective 
$n$-space extends stably to the projective $N$-space, $N > n$,    
if there exists a vector bundle on the larger space 
whose restriction to the smaller one is isomorphic to 
$E$ plus a direct sum of line bundles. We show that $E$ extends stably 
to the projective $N$-space for every $N > n$  
if and only if $E$ is the cohomology of a free monad (with three terms). 
The proof uses the method of Coand\u{a} and Trautmann (2006). 
Combining this result with a theorem of Mohan Kumar, Peterson and Rao 
(2003), we get a new effective version of the Babylonian tower 
theorem for vector bundles on projective spaces. 
\end{abstract}

\maketitle

Let $\mathbb{P}^n$, $n\geq 2$, be the projective $n$-space over an 
algebraically closed field $k$ of arbitrary characteristic, and let 
$S = k[X_0,\dots ,X_n]$ be the projective coordinate ring of 
$\mathbb{P}^n$. For an $m\geq 1$, let $R = k[X_0,\dots ,X_{n+m}]$ be the 
projective coordinate ring of $\mathbb{P}^{n+m}$ and embed $\mathbb{P}^n$ 
into $\mathbb{P}^{n+m}$ as the linear subspace $L$ of equations 
$X_{n+1} = \dots = X_{n+m} = 0$. One says that a vector bundle (= locally 
free coherent sheaf) $E$ on $\mathbb{P}^n$ {\it extends} to 
$\mathbb{P}^{n+m}$ if there exists a vector bundle $F$ on $\mathbb{P}^{n+m}$ 
such that $F\, \vert\, \mathbb{P}^n \simeq E$. The Babylonian tower 
theorem of Barth and Van de Ven \cite{bvv}, Tyurin \cite{tyu}, and  
E. Sato \cite{sa1}, \cite{sa2} (extended by Flenner \cite{fle} to vector 
bundles on the punctured spectrum of a regular local ring) asserts that a 
vector bundle on $\mathbb{P}^n$ can be extended to $\mathbb{P}^{n+m}$, 
$\forall m\geq 1$, if and only if it is a direct sum of line bundles. This 
result was (implicitly) conjectured by Horrocks \cite{ho2}. In order to 
attack this conjecture by using the techniques developed in \cite{ho1}, 
Horrocks introduced the following definition: $E$ {\it extends stably} to 
$\mathbb{P}^{n+m}$ if there exists a vector bundle $F$ on 
$\mathbb{P}^{n+m}$ and a (finite) direct sum of line bundles $A$ on 
$\mathbb{P}^n$ such that $F\, \vert \, \mathbb{P}^n \simeq E\oplus A$, and 
proved the following result:   

\vskip3mm
{\bf 0. Theorem.} (Horrocks \cite{ho2}, Theorem 1(ii))\quad 
\textit{If} $\text{H}^1_{\ast}(E) = 0$ \textit{and} 
$\text{H}^1_{\ast}(E^{\ast}) = 0$ \textit{and if} $E$ \textit{extends stably 
to} $\mathbb{P}^{2n-3}$ \textit{then} $E$ \textit{is a direct sum of line 
bundles}. 
\vskip3mm 

Recall that, for $0\leq i\leq n$, $\text{H}^i_{\ast}(E)$ denotes the graded 
$S$-module $\bigoplus_{d\in \mathbb{Z}}\text{H}^i(E(d))$. Notice, also, that, 
by Serre duality, the condition $\text{H}^1_{\ast}(E^{\ast}) = 0$ is 
equivalent to the condition $\text{H}^{n-1}_{\ast}(E) = 0$. 
  
The aim of this note is to characterize the vector bundles on 
$\mathbb{P}^n$ which are {\it infinitely stably extendable}. This means that, 
$\forall m\geq 1$, there exists a vector bundle $E^{(m)}$ on 
$\mathbb{P}^{n+m}$ and a direct sum of line bundles $A^{(m)}$ on 
$\mathbb{P}^n$ such that $E^{(m)}\, \vert \, \mathbb{P}^n \simeq 
E \oplus A^{(m)}$. An example of such a bundle is the cotangent bundle 
${\Omega}^1_{\mathbb{P}^n}$ of $\mathbb{P}^n$. The above mentioned 
characterisation is the content of the following: 

\vskip3mm 
{\bf 1. Theorem.}\quad 
\textit{For a vector bundle} $E$ \textit{on} $\mathbb{P}^n$, \textit{the 
following conditions are equivalent}:\\
\hspace*{3mm} (i) $\text{H}^i_{\ast}(E) = 0$ \textit{for} $1 < i < n-1$,\\
\hspace*{3mm} (ii) $E$ \textit{is the cohomology of a free monad on} 
$\mathbb{P}^n$,\\
\hspace*{3mm} (iii) $E$ \textit{is infinitely stably extendable}. 
\vskip3mm 

Recall that a {\it free monad} on $\mathbb{P}^n$ is a complex $K^{\bullet}$ 
consisting of (finite) direct sums of line bundles on $\mathbb{P}^n$, 
with $K^i = 0$ for $i < -1$ and for $i > 1$, with $d^{-1}_K : K^{-1} 
\rightarrow K^0$ a locally split monomorphism, and with $d^0_K : K^0 
\rightarrow K^1$ an epimorphism. The only non-zero cohomology sheaf 
$\mathcal{H}^0(K^{\bullet})$ of $K^{\bullet}$ is a vector bundle on 
$\mathbb{P}^n$. 

Using Theorem 1 and the main result of the paper of 
Mohan Kumar, Peterson and Rao \cite{mpr}, 
which asserts that a vector bundle satisfying condition (i) 
of Theorem 1 has rank $\geq n-1$ unless it is a direct sum of line bundles,  
we shall obtain, in Theorem 7 below, a 
new effective version of the Babylonian tower theorem. 

In order to prove Theorem 1, we shall use Horrocks' approach (which is 
very elementary) and the method developed by G. Trautmann and the author 
in \cite{ctr} for the proof of an effective version of the Babylonian 
tower theorem (which improves the estimate of Flenner \cite{fle}). 
Horrocks' proof of Theorem 0 is based on the following: 

\vskip3mm
{\bf 2. Lemma.}\quad 
\textit{Assume that the vector bundle} $E$ \textit{on} $\mathbb{P}^n$ 
\textit{extends to a vector bundle} $F$ \textit{on} $\mathbb{P}^{n+1}$. 
\textit{Then}, \textit{for a fixed} $1\leq i \leq n-1$, 
$\text{H}^i_{\ast}(E) = 0$ \textit{if and only if} 
$\text{H}^i_{\ast}(F) =0$ \textit{and} $\text{H}^{i+1}_{\ast}(F) = 0$. 
\vskip3mm 

For the proof of Lemma 2 one uses the short exact sequence (of sheaves on 
$\mathbb{P}^{n+1}$) :  $0\rightarrow F(-1)\rightarrow F\rightarrow E
\rightarrow 0$, the fact that $\text{H}^i(F(d)) = 0$ for $d << 0$ and 
$\text{H}^{i+1}(F(d)) = 0$ for $ d >> 0$. With the same kind of argument 
one can prove: 

\vskip3mm 
{\bf 3. Lemma.}\quad 
\textit{Under the hypothesis of Lemma} 2, $\text{H}^i_{\ast}(F)\rightarrow 
\text{H}^i_{\ast}(E)$ \textit{is surjective if and only if} 
$\text{H}^{i+1}_{\ast}(F) = 0$. 
\vskip3mm 

Horrocks' argument is reflected in the proof of the following: 

\vskip3mm 
{\bf 4. Lemma.}\quad 
\textit{Assume that the vector bundle} $E$ \textit{on} $\mathbb{P}^n$, 
$n\geq 4$, \textit{extends to a vector bundle} $F$ \textit{on} 
$\mathbb{P}^{n+m}$, \textit{for some} $m\geq n-3$. \textit{If there exists 
an} $(n+1)$-\textit{dimensional linear subspace} $P$ \textit{of} 
$\mathbb{P}^{n+m}$, \textit{containing} $L = \mathbb{P}^n$, 
\textit{such that the morphisms} $\text{H}^1_{\ast}(F\, \vert \, P) 
\rightarrow \text{H}^1_{\ast}(E)$ \textit{and} $\text{H}^1_{\ast}(F^{\ast}\, 
\vert \, P) \rightarrow \text{H}^1_{\ast}(E^{\ast})$ \textit{are surjective}, 
\textit{then} $\text{H}^i_{\ast}(E) = 0$ \textit{for} $1 < i < n-1$. 
\vskip3mm 

\begin{proof} 
Applying Lemma 3 to $E$ and $F\, \vert \, P$ one gets that $\text{H}^2_{\ast}
(F\, \vert \, P) =0$ and $\text{H}^2_{\ast}(F^{\ast}\, \vert \, P) = 0$ 
hence, by Serre duality, $\text{H}^{n-1}_{\ast}(F\, \vert \, P) = 0$. 
Considering a flag $P=P^1\subset P^2\subset \dots \subset P^m = 
\mathbb{P}^{n+m}$ of linear subspaces of $\mathbb{P}^{n+m}$ and applying 
Lemma 2 upstairs, one deduces that $\text{H}^i_{\ast}(F) = 0$, for 
$2 \leq i \leq n+m-2$. Considering, now, the flag $\mathbb{P}^{n+m} = P^m 
\supset \dots \supset P^1 = P\supset L = \mathbb{P}^n$ and applying Lemma 2 
downstairs, one obtains that $\text{H}^i_{\ast}(E) = 0$, for 
$2 \leq i \leq n-2$. 
\end{proof} 
\vskip3mm 

In order to prove the existence of a linear subspace $P$ satisfying the 
hypothesis of Lemma 4, we use the method developed in \cite{ctr}. 
Let $L^{\prime}$ be the linear subspace of $\mathbb{P}^{n+m}$ of 
equations $X_0 = \dots = X_n = 0$, with coordinate ring 
$S^{\prime} = k[X_{n+1},\dots ,X_{n+m}]$, and let $\pi : \mathbb{P}^{n+m} 
\setminus L^{\prime} \rightarrow L = \mathbb{P}^n$ be the linear projection. 
Let $L_i$ denote the $i$th infinitesimal neighbourhood of $L$ in 
$\mathbb{P}^{n+m}$, defined by the ideal sheaf ${\mathcal I}^{i+1}_L$. 
$\pi \, \vert \, L_i \rightarrow L$ is a retract of the inclusion 
$L \hookrightarrow L_i$ and endows ${\mathcal O}_{L_i}$ with the structure 
of an ${\mathcal O}_L$-algebra. As an ${\mathcal O}_L$-module: 
\[
{\mathcal O}_{L_i}\, \simeq \, {\mathcal O}_L\oplus 
{\mathcal O}_L(-1)\otimes_kS^{\prime}_1\oplus \dots \oplus 
{\mathcal O}_L(-i)\otimes_kS^{\prime}_i .
\]

Now, using the method from \cite{ctr}, we can prove the following result, 
which is the key point of the proof of Theorem 1: 

\vskip3mm 
{\bf 5. Lemma.}\quad 
\textit{Let} $E$ \textit{be a vector bundle on} $\mathbb{P}^n$, $n \geq 3$, 
\textit{and assume that} $E$ \textit{can be extended to a vector bundle} 
$F$ \textit{on} $\mathbb{P}^{n+m}$, \textit{for some} $m \geq 1$. 
\textit{Let} $\xi$ \textit{be a fixed element of} $\text{H}^1(E)$. 
\textit{Then there exists a homogeneous ideal} $J \subseteq S^{\prime}_+$, 
\textit{with at most} $\text{h}^2(E(-i)) := \text{dim}_k\text{H}^2(E(-i))$ 
\textit{minimal homogeneous generators in degree} $i$, $\forall i \geq 1$, 
\textit{such that if} $p$ \textit{is a point of the subscheme of} 
$L^{\prime}$ \textit{defined by} $J$, \textit{and if} $P$ 
\textit{is the linear subspace of} $\mathbb{P}^{n+m}$ \textit{spanned by} 
$L$ \textit{and} $p$, \textit{then} $\xi$ \textit{can be lifted to} 
$\text{H}^1(F\, \vert \, P)$. 
\vskip3mm 

\begin{proof} 
We shall construct recursively the homogeneous components $J_0 = (0), J_1, 
\dots , J_i, \dots $ of $J$. Assume that we have constructed $J_0 = (0), 
\dots ,J_i$ such that, if $Y_i$ is the closed subscheme of $L_i$ defined by 
the ideal sheaf: 
\[
{\mathcal O}_L(-1)\otimes_kJ_1\oplus \dots \oplus {\mathcal O}_L(-i)
\otimes_kJ_i, 
\]
$\xi \in \text{H}^1(E)$ lifts to an element ${\xi}_i\in 
\text{H}^1(F\, \vert \, Y_i)$. Let $Y^{\prime}_{i+1}$ be the closed subscheme 
of $L_{i+1}$ defined by the ideal sheaf: 
\[
{\mathcal O}_L(-1)\otimes_kJ_1\oplus \dots \oplus {\mathcal O}_L(-i)
\otimes_kJ_i\oplus {\mathcal O}_L(-i-1)\otimes_k(S^{\prime}_1J_i).
\]
Tensorizing by $F$ the exact sequence: 
\[
0\rightarrow {\mathcal O}_L(-i-1)\otimes_k(S^{\prime}_{i+1}/S^{\prime}_1J_i) 
\longrightarrow {\mathcal O}_{Y^{\prime}_{i+1}} \longrightarrow 
{\mathcal O}_{Y_i} \rightarrow 0
\]
and taking cohomology, one deduces an exact sequence: 
\[
\text{H}^1(F\, \vert \, Y^{\prime}_{i+1})\rightarrow 
\text{H}^1(F\, \vert \, Y_i)\overset{{\delta}^{\prime}}{\longrightarrow} 
\text{H}^2(E(-i-1))\otimes_k(S^{\prime}_{i+1}/S^{\prime}_1J_i).
\]
Choosing a basis of the $k$-vector space $\text{H}^2(E(-i-1))$, one can find 
a $k$-vector space $J_{i+1}$, with $S^{\prime}_1J_i \subseteq J_{i+1} 
\subseteq S^{\prime}_{i+1}$, and with $\text{dim}_k(J_{i+1}/S^{\prime}_1J_i) 
\leq \text{h}^2(E(-i-1))$, such that ${\delta}^{\prime}({\xi}_i) \in 
\text{H}^2(E(-i-1)) \otimes_k (J_{i+1}/S^{\prime}_1J_i)$. Let $Y_{i+1}$ be 
the closed subscheme of $L_{i+1}$ defined by the ideal sheaf: 
\[
{\mathcal O}_L(-1)\otimes_kJ_1\oplus \dots \oplus {\mathcal O}_L(-i)
\otimes_kJ_i\oplus {\mathcal O}_L(-i-1)\otimes_kJ_{i+1}.
\]
One has $Y_i\subseteq Y_{i+1}\subseteq Y^{\prime}_{i+1}$, and from the 
commutativity of the diagram: 
\[
\begin{CD}
\text{H}^1(F\, \vert \, Y^{\prime}_{i+1}) @>>> 
\text{H}^1(F\, \vert \, Y_i) @>{{\delta}^{\prime}}>> 
\text{H}^2(E(-i-1))\otimes_k(S^{\prime}_{i+1}/S^{\prime}_1J_i)\\
   @VVV @V{\text{id}}VV @VVV\\
\text{H}^1(F\, \vert \, Y_{i+1}) @>>> 
\text{H}^1(F\, \vert \, Y_i) @>{\delta}>> 
\text{H}^2(E(-i-1))\otimes_k(S^{\prime}_{i+1}/J_{i+1})
\end{CD}
\]
one deduces that $\delta ({\xi}_i) = 0$, hence ${\xi}_i$ lifts to an 
element ${\xi}_{i+1} \in \text{H}^1(F\, \vert \, Y_{i+1})$. 

Finally, assume that the above construction of $J$  has been completed and 
let $p$ be a point of the closed subscheme of $L^{\prime}$ defined by 
$J$. If $J(p) \subset S^{\prime}_+$ is the homogeneous ideal of the point 
$p$ of $L^{\prime} \simeq \mathbb{P}^{m-1}$, then the closed subscheme of 
$L_i$ defined by the ideal sheaf:
\[
{\mathcal O}_L(-1)\otimes_kJ(p)_1\oplus \dots \oplus 
{\mathcal O}_L(-i)\otimes_kJ(p)_i
\]
is exactly $P\cap L_i$, i.e., the $i$th infinitesimal neighbourhood in 
$P$ of the hyperplane $L$ of $P \simeq \mathbb{P}^{n+1}$. Since 
$J \subseteq J(p)$, it follows that $P \cap L_i \subseteq Y_i$, hence 
$\xi \in \text{H}^1(E)$ lifts to the image ${\xi}_i(p)$ of 
${\xi}_i \in \text{H}^1(F\, \vert \, Y_i)$ into 
$\text{H}^1(F\, \vert \, P\cap L_i)$. Using the exact sequence: 
\[
\text{H}^1(F\, \vert \, P)\longrightarrow 
\text{H}^1(F\, \vert \, P\cap L_i)\longrightarrow 
\text{H}^2((F\, \vert \, P)(-i-1))
\]
and the fact that $\text{H}^2((F\, \vert \, P)(-i-1)) = 0$ for 
$i >> 0$, one deduces that $\text{H}^1(F\, \vert \, P)\rightarrow 
\text{H}^1(F\, \vert \, P\cap L_i)$ is surjective for $i >> 0$, hence 
$\xi$ lifts to an element of $\text{H}^1(F\, \vert \, P)$.
\end{proof} 
\vskip3mm 

Before proving Theorem 1, we need one more remark, related to Lemma 4: 

\vskip3mm 
{\bf 6. Lemma.}\quad 
\textit{Under the hypothesis of Lemma} 4, \textit{the condition} 
``$\text{H}^1_{\ast}(F\, \vert \, P)\rightarrow \text{H}^1_{\ast}(E)$ 
\textit{surjective}'' \textit{is an open condition on the variety 
parametrizing the} $(n+1)$-\textit{dimensional linear subspaces} $P$ 
\textit{of} $\mathbb{P}^{n+m}$ \textit{containing} $L = \mathbb{P}^n$. 
\vskip3mm 

\begin{proof} 
Let $\widetilde{\mathbb{P}}$ be the subvariety of $\mathbb{P}^{n+m} \times 
L^{\prime}$ consisting of the pairs $(x,p)$ such that $x$ belongs to the 
linear span of $L$ and $p$. Actually, $\widetilde{\mathbb{P}}$ is the 
blow-up of $L$ in $\mathbb{P}^{n+m}$. Consider the incidence diagram: 
\[
\begin{CD}
\mathbb{P}^{n+m} @<{\sigma}<< {\widetilde{\mathbb{P}}} 
@>{\widetilde \pi}>> L^{\prime}.
\end{CD}
\]
$\widetilde{\mathbb{P}}$ is a $\mathbb{P}^{n+1}$-bundle over $L^{\prime}$ 
and, $\forall p\in L^{\prime}$, $\sigma$ maps isomorphically 
${\widetilde \pi}^{-1}(p)$ onto the linear span of $L$ and $p$. In this 
way, $L^{\prime}$ parametrizes the linear subspaces $P$ from the statement. 

Now, by Lemma 3, $\text{H}^1_{\ast}(F\, \vert \, P)\rightarrow 
\text{H}^1_{\ast}(E)$ is surjective if and only if 
$\text{H}^2_{\ast}(F\, \vert \, P) = 0$. There exists integers $d_1$, $d_2$ 
such that $\text{H}^1(E(d)) = 0$ for $d > d_1$ and $\text{H}^2(E(d)) = 0$ 
for $d \leq d_2$. Then, for {\it every} $(n+1)$-dimensional linear subspace 
$P$ of $\mathbb{P}^{n+m}$ containing $L = \mathbb{P}^n$, 
$\text{H}^2((F\, \vert \, P)(d)) = 0$ for $d \geq d_1$ and for 
$d \leq d_2$. Finally, if there exists a $P_0$ such that 
$\text{H}^2_{\ast}(F\, \vert \, P_0) = 0$ then a semi-continuity argument 
applied to ${\sigma}^{\ast}F$ over $L^{\prime}$ shows that, for $P$ in a  
neighbourhood of $P_0$, $\text{H}^2((F\, \vert \, P)(d)) = 0$ for 
$d_2 < d < d_1$, hence $\text{H}^2_{\ast}(F\, \vert \, P) = 0$. 
\end{proof} 
\vskip3mm 

\begin{proof}[Proof of Theorem 1] 
(i)$\Rightarrow$(ii) This is the well known method of Horrocks \cite{ho3}, 
\cite{ho4} of ``killing cohomology'' (see Barth and Hulek \cite{bah} for a 
more accessible reference) combined with the splitting criterion of 
Horrocks. 

(ii)$\Rightarrow$(iii) Let $K^{-1}\overset{d^{-1}}{\longrightarrow} K^0
\overset{d^0}{\longrightarrow} K^1$ be a free monad on $\mathbb{P}^n$. 
Embed $\mathbb{P}^n$ in $\mathbb{P}^{n+1}$ as the hyperplane of equation 
$X_{n+1} = 0$, let $p$ be the point $[0:\ldots :0:1]$ of $\mathbb{P}^{n+1}$ 
and let $\pi : \mathbb{P}^{n+1}\setminus \{p\} \rightarrow \mathbb{P}^n$ be 
the linear projection. ${\pi}^{\ast}K^i$ extends to a direct sum of line 
bundles ${\widetilde K}^i$ on $\mathbb{P}^{n+1}$ and ${\pi}^{\ast}d^i$ 
extends to a morphism ${\widetilde d}^{\, i} : {\widetilde K}^i \rightarrow 
{\widetilde K}^{i+1}$. One can consider on $\mathbb{P}^{n+1}$ the free 
monad $K^{\prime -1}\overset{d^{\prime -1}}{\longrightarrow} K^{\prime 0} 
\overset{d^{\prime 0}}{\longrightarrow} K^{\prime 1}$, where: 
\begin{gather*} 
K^{\prime -1}={\widetilde K}^{-1},\  \  
K^{\prime 0}={\widetilde K}^{-1}(1)\oplus {\widetilde K}^0\oplus 
{\widetilde K}^1(-1),\  \  K^{\prime 1}={\widetilde K}^1\\
d^{\prime -1}=(X_{n+1}\cdot -\  ,\   {\widetilde d}^{\, -1}\  ,\   
0)^{\text{trsp}}, \  \  
d^{\prime 0}=(0\  ,\  {\widetilde d}^{\, 0}\  ,\  X_{n+1}\cdot -).
\end{gather*}
Then: ${\mathcal H}^0(K^{\prime \bullet})\, \vert \, \mathbb{P}^n \simeq 
K^{-1}(1)\oplus {\mathcal H}^0(K^{\bullet})\oplus K^1(-1)$. 

(iii)$\Rightarrow$(i) Assume that a vector bundle $E$ on ${\mathbb P}^n$, 
$n \geq 4$, extends stably to a vector bundle $F$ on ${\mathbb P}^{n+m}$, 
for some $m$ large enough (this will be made precise below). We want to 
show that, in this case, $\text{H}^i_{\ast}(E) = 0$ for $1 < i < n-1$. 
By definition, $F\, \vert \, {\mathbb P}^n \simeq E\oplus A$, where $A$ 
is a direct sum of line bundles. Since $\text{H}^i_{\ast}(A) = 0$ for 
$1 \leq i \leq n-1$, it follows that $\text{H}^i_{\ast}(E\oplus A) 
\simeq \text{H}^i_{\ast}(E)$, for $1 \leq i \leq n-1$. 

Let ${\mu}_i := \text{dim}_k(\text{H}^1(E(i))/S_1\text{H}^1(E(i-1)))$ be 
the number of minimal homogeneous generators of degree $i$ of the graded 
$S$-module $\text{H}^1_{\ast}(E)$. Choosing a minimal system of homogeneous 
generators of the graded $S$-module $\text{H}^1_{\ast}(E\oplus A)$ and 
applying Lemma 5 to each of these generators, one deduces that, if 
$m > \sum_{i>j}{\mu}_i\text{h}^2(E(j))$, there exists an $(n+1)$-dimensional 
linear subspace $P_0$ of ${\mathbb P}^{n+m}$, containing $L = 
{\mathbb P}^n$, such that $\text{H}^1_{\ast}(F\, \vert \, P_0) 
\rightarrow \text{H}^1_{\ast}(E\oplus A)$ is surjective. It follows now, 
from Lemma 6, that for the {\it general} $(n+1)$-dimensional linear 
subspace $P$ of ${\mathbb P}^{n+m}$, containing $L = {\mathbb P}^n$, 
$\text{H}^1_{\ast}(F\, \vert \, P) \rightarrow \text{H}^1_{\ast}
(E\oplus A)$ is surjective. 

Similarly, putting ${\mu}_i^{\ast} := \text{dim}_k(\text{H}^1(E^{\ast}(i))
/S_1\text{H}^1(E^{\ast}(i-1)))$, it follows that, for $m > \sum_{i>j}
{\mu}_i^{\ast}\text{h}^2(E^{\ast}(j))$, for the {\it general} $P$ as above, 
$\text{H}^1_{\ast}(F\, \vert \, P) \rightarrow \text{H}^1_{\ast}
(E^{\ast}\oplus A^{\ast})$ is surjective. 

One deduces now, from Lemma 4, that if: 
\[
m > \text{max}(\underset{i>j}{\textstyle \sum} {\mu}_i\text{h}^2(E(j)),\  
\underset{i>j}{\textstyle \sum} {\mu}_i^{\ast}\text{h}^2(E^{\ast}(j)),\  
n-4) 
\]
then $\text{H}^i_{\ast}(E\oplus A) = 0$ for $1 < i < n-1$, hence 
$\text{H}^i_{\ast}(E) = 0$, for $1 < i < n-1$.    
\end{proof}
\vskip3mm
 
Using the estimate from the last part of the proof of Theorem 1 and the 
main result of the paper of Mohan Kumar et al. \cite{mpr}, one gets the 
following version of the Babylonian tower theorem (compare with \cite{bvv}, 
\cite{tyu}, \cite{sa1}, \cite{sa2}, \cite{fle}, and \cite{ctr}): 

\vskip3mm 
{\bf 7. Theorem.}\quad 
\textit{Let} $E$ \textit{be a vector bundle of rank} $r$ \textit{on} 
${\mathbb P}^n$, $n \geq 4$. \textit{If} $E$ \textit{extends to a vector 
bundle} $F$ \textit{on} ${\mathbb P}^{n+m}$, \textit{for some} 
\[
m > \text{max}(\underset{i>j}{\textstyle \sum} {\mu}_i\text{h}^2(E(j)),\  
\underset{i>j}{\textstyle \sum} {\mu}_i^{\ast}\text{h}^2(E^{\ast}(j)),\  
n-4,\  r-n+1), 
\]
(\textit{see the proof of Theorem} 1 \textit{for the notation}) 
\textit{then} $E$ \textit{splits as a direct sum of line bundles}. 
\vskip3mm 

\begin{proof} 
By (the proof of) Theorem 1, $\text{H}^i_{\ast}(E) = 0$ for $1 < i < n-1$, 
hence, by Lemma 2, $\text{H}^i_{\ast}(F) = 0$, for $1 < i < n+m-1$. Since 
$\text{rk}F = r < n+m-1$, the main result of \cite{mpr} implies that 
$F$ is a direct sum of line bundles, hence the same is true for 
$E \simeq F\, \vert \, {\mathbb P}^n$. 
\end{proof}

\end{document}